\documentclass[10pt,reqno]{amsart}

\usepackage[hypertexnames=false, colorlinks, citecolor=red, linkcolor=red ]{hyperref}
\usepackage{geometry,amsmath,amssymb,enumerate,bbm}
\usepackage[american]{babel}
\newcommand{\assign}{:=}

\newcommand{\tmem}[1]{{\em #1\/}}
\newcommand{\tmop}[1]{\ensuremath{\operatorname{#1}}}

\newenvironment{enumeratealpha}{\begin{enumerate}[a{\textup{)}}] }{\end{enumerate}}

\usepackage{graphicx}

\usepackage{psfrag}

\usepackage{amsmath} 

\usepackage{amsfonts}

\usepackage{amssymb}

\usepackage{amsthm} 

\usepackage{tikz}
\usepackage{tikz-qtree}
\usepackage{float}
\newtheorem{thm}{Theorem}[section]
\newtheorem{prop}[thm]{Proposition}

\newtheorem{lemma}[thm]{Lemma}
\newtheorem{defn}[thm]{Definition}

\newtheorem{theorem}{Theorem}

\newcommand{\HHH}{{\mathcal{H}}}
\newcommand{\III}{{\mathcal{I}}}

\newcommand{\diam}{\text{diam}}

\newcommand{\Stop}{\text{Stop}}
\newcommand{\icirc}{{~^o\!I}}

\begin{document}
\title{Ahlfors regular spaces have regular \\subspaces of any dimension}

\author[N. Arcozzi, A. Monguzzi, M. Salvatori]{Nicola Arcozzi*, Alessandro Monguzzi**, Maura Salvatori***}

\address{*Dipartimento di Matematica, Alma Mater Studorium Universit\`a di
  Bologna, Piazza di Porta San Donato 5, 40126 Bologna, Italy}
\email{{\tt nicola.arcozzi@unibo.it}}

 \address{**Dipartimento di Matematica e Applicazioni, Universit\`a degli Studi di Milano--Bicocca, Via R. Cozzi 55, 20126 Milano, Italy}  
\email{{\tt alessandro.monguzzi@unimib.it}}

\address{***Dipartimento di Matematica, Universit\`a degli Studi di
  Milano, Via C. Saldini 50, 20133 Milano, Italy}

\email{{\tt maura.salvatori@unimi.it}}

\thanks{All authors are members of INDAM-GNAMPA. The first and second authors are partially supported by the grant GNAMPA 2020 \emph{Alla frontiera tra l'analisi complessa in più variabili e l'analisi armonica}. The third author is partially supported by the grant GNAMPA 2020 \emph{Fractional Laplacians and subLaplacians on Lie groups and trees}. The first author is partially supported by
European Unions Horizon 2020 research programme, Marie Sk\l{}odowska-Curie grant agreement No. 777822}


\date{\today}

\begin{abstract}
We characterize $Q$-dimensional Ahlfors regular spaces among trees' boundaries and show how to construct, for each $0 < \alpha < Q$, an $\alpha$-regular subspace. As an application, we give an alternative simple proof of the existence of $\alpha$-regular subspaces of a $Q$-dimensional complete Ahlfors regular metric space $(X,\rho)$, which was proved in \cite{JJKRRS}.
\end{abstract}
\maketitle

Any metric space $(X,\rho)$ supports a $Q$-dimensional Hausdorff measure $\HHH^Q$ for all $Q>0$, see \cite{Rogers}
for a nice introduction to the topic. The upper Hausdorff dimension of $X$ is 
$Q=\text{dim}(X):=\inf\{\alpha:\ \HHH^\alpha(X)=\infty\}\ge0$. If $0<\alpha<Q$, it is natural to ask whether there 
exist subsets 
$A$ of $X$  such that $0<\HHH^\alpha(A)<\infty$. 
A special case of Corollary 7 in \cite{How95} gives a positive answer in great generality: 
if $X$ is a complete, separable Hausdorff space, and 
$0<\alpha<\text{dim}(X)$, then there is a compact subset of $A$ such that $0<\HHH^\alpha(A)<\infty$.

Among metric spaces, a class which has been much studied in recent decades is that of the 
\it Ahlfors regular \rm ones. In these spaces the $Q$-dimensional Hausdorff measure 
$\HHH^Q(B(x,r))$ of the metric ball having center at $x$ and radius $r>0$ satisfies
\begin{equation}\label{ahlfors}
cr^Q\le\HHH^Q(B(x,r))\le Cr^Q.
\end{equation}
with $0<c<C$ independent of $x$ and $r\le\text{diam}(X)\le\infty$. 
They proved to be an excellent setting for 
developing harmonic analysis even in many of its fine features \cite{DS}. In this paper we always assume $(X,\rho)$ to be also complete. 

Recently, Ahlfors regular spaces were found to be a favorable environment for many chapters of potential theory 
\cite{ARSW}. Adding more structure, rectifiable curves and ``gradients'' of functions, it is possible to 
push the potential theory even further \cite{HK}, but in this article we will make use of Ahlfors regularity alone.

It is thus desirable having a vast catalog of Ahlfors regular spaces. In this article we increase the catalog by characterizing all Ahlfors regular spaces among trees boundaries. Then, for each given complete Ahlfors regular tree boundary we construct ``fractals'' of any lower dimension, which are themselves complete Ahlfors regular. 
A different construction can be obtained by specializing Corollary 5.2 of \cite{JJKRRS} and the constructions leading to its proof to our context. In fact, the result in \cite{JJKRRS} can be obtained by its special ``dyadic'' version and Michael Christ's Theorem on the dyadicization of homogeneous spaces \cite{Christ}.

\begin{theorem}\label{theomain}
Let $( X, \rho)$ be a complete Ahlfors $Q$-regular metric space
and let $0 < \alpha < Q$. Then, there exists $Y \subseteq X$, closed, such that $( Y,
\rho)$ is $\alpha$-regular.
\end{theorem}
We point out that the result is classical if $(X,\rho)$ has a dilation structure. 
See e.g. \S8.3 in \cite{Fal}, from which a statement like ours can can be 
easily deduced.
In fact, the easiest way to produce fractal subsets of Euclidean space is by means of families of similarities. 
When available, dilations are the easiest 
way to produce a multi-scale decomposition of the metric space; but such decompositions do not exist in all Ahlfors 
regular space, and the basic idea 
is using multiscale decompositions to produce Ahlfors regular fractal sets of all dimensions.

Our approach is to split the proof into two steps. 
First, we show that Theorem \ref{theomain} holds if $X$ is a tree boundary 
(equivalently, if $X$ is an ultrametric space: $\rho(x,y)\le\max\{\rho(x,z),\rho(z,y)\}$). 
  This is done in Section \ref{alberi}. 
  While doing so, we will provide several characterizations of trees 
  whose boundary is Ahlfors regular, which might have independent interest. The general case is reduced to the special one 
  by means of  standard dyadization in Section \ref{dyadization}.  
  In fact, we see here at work the heuristic principle that statements concerning Ahlfors regular spaces in general 
  can be deduced by the corresponding statements restricted to regular tree boundaries, or regular ultrametric spaces.
  
  \vskip0.5cm
  
  \noindent \bf Notation. \rm We write $\sharp A$ to denote the cardinality of $A$. 
  We use the symbols $\lesssim$ and $\approx$ with their usual meaning: $A(x)\lesssim B(x)$ 
  means that $A(x)\le C B(x)$ with a constant 
  $C>0$ which is independent of $x$ (but which might depend on other parameters); $A(x)\approx B(x)$ means 
  $A(x)\lesssim B(x)\lesssim A(x)$.

  \section{Ahlfors regular tree boundaries.}\label{alberi}
  
  A \it tree \rm $T$  is a connected, simply connected graph, with edge set $E(T)$ and vertex set $V(T)$, which, 
  with harmless abuse, we 
  denote by $T$. The edge 
having as endpoints the vertices $x\ne y$ can be identified with the unordered couple $\{x,y\}$. A path between two
vertices $x$ and $y$ is a sequence of edges $\{x_0=x,x_1\},\ \{x_1,x_2\}, \dots, \{x_n-1,x_n=y\}$: $n\ge0$ is the 
length of the path.
The natural distance $|y-x|$ between two vertices $x,y$ in $T$ is the minimal length of a path starting at $x$ and 
ending at $y$. It is realized by
exactly one path $[x,y]$, the geodesic between $x$ and $y$. It is correct and useful to think of $|y-x|$ as the 
hyperbolic distance between $x$ and $y$.

We now introduce a well known analog of the Euclidean metric, sometimes called the visual metric.
A distinguished vertex $o$ is chosen: the \textit{root} of $T$.  We write $y\succ x$ if $x\in[o,y]$, and 
 $|x|:=|x-o|$. Given vertices $x,y$, we 
have that $[o,x]\cap[o,y]=[o,x\wedge y]$
for exactly one vertex $x\wedge y$, the \it confluent \rm of $x$ and $y$. If $e=\{x,y\}$ is an edge of $T$ and $y$ 
is its endpoint which is furthest
from $o$, set ${\bar\rho}(e)=2^{-|y|-1}$. Given a finite path $\Gamma=\{e_j\}_{j=1}^n$ in $T$, define its ${\bar\rho}$-length to be 
${\bar\rho}(\Gamma)=\sum_{j=1}^n {\bar\rho}(e_j)$.
The ${\bar\rho}$-distance ${\bar\rho}(x,y)$ between the vertices $x$ and $y$ is the minimal ${\bar\rho}$-length of a path starting at 
$x$ and ending at $y$. It is obviously
realized by the ${\bar\rho}$-length of the geodesic $[x,y]$. The metric space $(T,{\bar\rho})$ has finite diameter. 

We shall make on $T$ the assumption that the number of 
the edges leaving a vertex $x\ne o$ is at least $2$. This rules out the existence of ``leaves'' 
at a finite distance from the root. We will see that $Q$-regularity of $T$'s boundary implies that the number of such 
edges is bounded by a constant $N$ which is independent of the particular vertex..
Let $\overline{T}\supset T$ 
be the completion of $T$ with respect to the metric ${\bar\rho}$: $(\overline{T},{\bar\rho})$ is a complete metric space; 
$(T,{\bar\rho})$ is dense into it and it
is a discrete set. We define $\partial T:=\overline{T}\setminus T$ the \it boundary \rm of $T$. The points $\zeta$ in 
$\partial T$ are in a bijection with half-infinite geodesics 
$\Gamma$ starting at $o$ (i.e. infinite paths starting at $o$, with no repeated edges), since the sequence of the 
vertices in $\Gamma$ is a Cauchy sequence
in $(T,{\bar\rho})$, converging to a unique $\zeta$ in $\partial T$. In this case we write $\Gamma_\zeta=[o,\zeta)$ 
and $\overline{\Gamma}_\zeta=[o,\zeta]$.

Given a tree $T$ with a root $o$, the subtree $T(x)$ having root $x$ has as vertices the $y$'s 
such that $[o,y]\ni x$. Its boundary, $\partial T(x)$, can be naturally identified with a subset of $\partial T$.

The restriction of the distance $\bar\rho$ to $\partial T=\bar T\setminus T$ is
\[ \bar\rho ( x, y) = 2^{- | x \wedge y |}, \]
and an inspection shows that is an ultrametric, $\rho(x,y)\le\min\{\rho(x,z),\rho(z,y)\}$ on $\partial T$. 
In particular, any point of a metric ball $B$ in $\partial T$
is its center.
If the root has at least two descendants, 
$\diam_{\bar\rho} ( \partial T) = 1$. Here  $\wedge$ is the extension of the confluent to $\partial T$: 
$[o,x\wedge y]=\Gamma_x\cap\Gamma_y$.
Balls $\overline{B_{\bar\rho} ( x, 2^{- n})} = B_{\bar\rho} \left( x, 2^{- n +
\frac{1}{2}} \right)$ in the metric are cl-open. Metric balls can be identified with the boundaries of the subtrees 
$\partial T(w)$:
$$
\overline{B_{\bar\rho} ( x, 2^{- n})} = B_{\bar\rho} \left( x, 2^{- n +\frac{1}{2}} \right)=\partial T(w)
$$
for a well specified $w$ in $T$ with $|w|=n$. We are not ruling out the possibility that 
$\overline{B_{\bar\rho} ( x, 2^{- n})} =\overline{B_{\bar\rho} ( x, 2^{- n-1})} $, although this, in the $Q$-regular case, 
can not happen too often, as can be made quantitative using the Lemma below. 

We provide here two characterizations of trees whose boundaries 
are $Q$-regular. 

Given $x\in T$, we denote by $\tmop{Stop}(x)$ the family of stopping regions below $x$, that is,  
$T\supseteq S\in\tmop{Stop} (x)$ if 
\begin{equation}\label{stopping-region}
S\subset T(x) \textrm{ and } \bigsqcup_{y\in S}\partial T(y)=\partial T(x).
\end{equation}
where  $\bigsqcup$ denotes a disjoint union. 
Stopping regions $S$ parametrize open covers of $\partial T(x)$ by metric balls. Since $\partial T(x)$ is compact, $S$ 
is \it a fortiori \rm finite.

\begin{lemma}\label{lemma-regular-tree}
  $( \partial T, \bar \rho)$ is $Q -$regular if and only if there are $0<a_0<a_1$ so that, for all $x$ in $T$ and $S$ in 
  $\tmop{Stop} ( x)$, we have
  \begin{equation}\label{lemma-reg-1} 
  a_0 2^{- Q | x |}\le\sum_{y \in S} 2^{- Q | y |} \le a_1 2^{- Q | x |} .
  \end{equation}
\end{lemma}
Applied to the stopping set $S(x)$ containing the vertices just below $x$, the lemma implies that $\sharp(S(x))\le N$, where 
$N$ is independent of $x$.

\begin{proof} Assume that $\partial T$ is $Q$-regular. We first show that if $\zeta\in\partial T$ and  
$\overline{B}(\zeta,2^{-n})=\overline{B}(\zeta,2^{-n-1})=\dots \overline{B}(\zeta,2^{-n-m})$, then $m\le C$. 
There are $x_j$ in $T$ with $\overline{B}(\zeta,2^{-n-j})=\partial T(x_j)$ and $|x_j|=n+j$.
By $Q$-regularity, then,
$$
c2^{-Qn}\le\HHH^Q(\partial T(x_0))=\HHH^Q(\partial T(x_m))\le C2^{-Q(n+m)}, 
$$
so that $2^{-Qm}\ge\frac{c}{C}$. We have then the estimate:
\begin{equation}\label{forgotten}
 \frac{c^2}{C}2^{-Q|x|}\le\HHH^Q(\partial T(x))\le C2^{-Q|x|}
\end{equation}
Let now $x\in T$ and $S\in\tmop{Stop}(x)$ be given and let 
$z \in \partial T(x) \subseteq \partial T$.
 Then, $\partial T ( z) = \sqcup_{y \in S} \partial T( y)$, 
hence,
\[
\frac{c^2}{C} 2^{- Q | x |} \le \mathcal{H}^Q ( \partial T ( x)) = \sum_{y \in S}
   \mathcal{H}^Q ( \partial T ( y)) \le C\sum_{y\in S}  2^{- Q | y |}, 
   \]
  where both the first and the last inequalities hold by $Q$-regularity \eqref{ahlfors},
thus 
$$
2^{- Q | x |}\le\frac{C^2}{c^2}\sum_{y\in S}  2^{- Q | y |}, 
\text{ and similarly } 2^{- Q | x |}\ge\frac{c^2}{C^2}\sum_{y\in S}  2^{- Q | y |}.
$$
  Conversely, let $M>0$ be a fixed integer and let $\{U_n\}_n$ be a cover of $\partial T(x)$ by sets $U_n\subseteq \partial T$ 
  such that $\diam(U_n)\leq 2^{-M}$. For each $n$ select $w_n \in U_n$ and replace $U_n$ by $\partial T(w_n)\supseteq U_n$ with 
  $2^{-|w_n|-1}<\diam(U_n)\le 2^{-|w_n|}\le 2^{-M}$. 
  Then, $\{\partial T(w_n)\}$ is a new cover of $\partial T(x)$ and, since $\partial T(x)$ is compact, 
  we can select a finite subfamily of it, which we might relabel as 
  $\{\partial T(w_n)\}_{n=1}^N$, such that $\{w_n\}$ is a stopping set below $x$. 

    Now, 
  \begin{align*}
   2^Q\sum_{n=1}^N \diam(U_n)^Q&\ge \sum_{n=1}^N \diam(\partial T(w_n))^Q\\
   &= \sum_{n=1}^N2^{-Q|w_n|}\\
   &\ge a_0\,2^{-Q|x|},
  \end{align*}
  and similarly
  $$
  \sum_{n=1}^N \diam(U_n)^Q\le a_1 2^{-Q|x|}.
  $$
In particular, given $\partial T(x)= \overline{B_{\bar\rho}(\zeta, 2^{-|z|})}$\,, we deduce that
\begin{eqnarray*}
a_0 2^{-Q} 2^{-Q|x|}&\le& \mathcal{H}^Q_{\delta} ( \partial T(x))  
 \assign  \inf \bigg\{ \sum_n \diam ( U_n)^Q : \cup_{n} U_n\supseteq \partial T(x),\, \tmop{diam} ( U_n) 
 \leqslant \delta \bigg\}\crcr
 &\le& a_1 2^{-Q|x|},
\end{eqnarray*}
  hence,
  $$
2^{-Q}a_0 2^{-Q|x|}\le  \mathcal H^{Q}(\partial T(x))=\lim_{\delta\to 0}\mathcal H^{Q}_\delta(\partial T(x))\le a_1 2^{-Q|x|}.
  $$
\end{proof}

Let $k\geq 1$ be an integer and, for any $x\in T$ set
$$
S_k(x):=\{y\in T(x): |y-x|=k\}.
$$
\begin{thm}\label{ahl-regularity-trees}
  The metric space $(\partial T,\bar\rho)$ is $Q$-regular if and only if there are $0<b_0<b_1$ so that
  \begin{equation}\label{combinatorial}
   b_0 2^{Q k}\le \sharp S_k ( x)\le b_1 2^{Q k}
  \end{equation}
  holds for all $x\in T$.
\end{thm}

\begin{proof} Suppose \eqref{combinatorial} holds.
By natural rescaling, it suffices to prove the result for $x=o$, the root of the tree $T$. Let
$S \in \tmop{Stop} ( o)$ be a stopping set, let $x \in S$  and consider $S_{N-|x|}(x)$ where $N=\max\{|y|:\ y\in S\}$. 
is a positive integer to be fixed later. 
Then, 
\begin{align*}
 2^{-QN}\sharp S_N(o)&=&\sum_{u\in S_N(o)}2^{-Q|u|}=\sum_{x\in S}\sum_{u\in S_{N-|x|}(x)}2^{-Q(|x|+N-|x|)}\crcr
 &=&\sum_{x\in S}2^{-Q|x|}\sharp S_{N-|x|}(x)2^{-Q(N-|x|)}.
\end{align*}
By \eqref{combinatorial} we have 
$$b_0\le 2^{-QN}\sharp S_N(o)\le b_1\sum_{x\in S}2^{-Q|x|}\le\frac{b_1}{b_0}2^{-QN}\sharp S_N(o)\le\frac{b_1^2}{b_0},$$
and we can apply Lemma \ref{lemma-regular-tree}.

Conversely, assume that $(\partial T,\bar\rho)$ is $Q$-regular. 
Then, since each $S_{k}(x)\in \Stop(x)$, Lemma \ref{lemma-regular-tree} implies that
$$
2^{-Q|x|}\approx \sum_{y\in S_k(x)}2^{-Q|y|}= 2^{-Q(|x|+k)}\sharp S_k(x),
$$
and the desired estimate for $\sharp S_k ( x)$ follows.
\end{proof}

By Theorem \ref{ahl-regularity-trees} the boundaries of homogeneous trees are Ahlfors regular 
metric spaces. Namely, if $T$ is a $q$-homogeneous tree, than $(\partial T,\bar\rho)$ is  $\log_2 q$-regular.

\vskip0.5cm

Before proving our main theorem for a general tree, it is simple and useful to see how it is proved 
in the model example of the binary tree. We exploit Theorem \ref{ahl-regularity-trees} above.

Let $T$ be the binary tree with root $o$ so that $(\partial T,\bar\rho)$ is $1$-regular. Given $\alpha\in ( 0, 1)$, we want to 
find a $\alpha$-regular subspace of 
$\partial T$. We construct a subtree $U$ of $T$, where a {\tmem{subtree}} of $T$ consists in a subset $E ( U)$ of the edges 
$E ( T)$ of
$T$ such that, together with their endpoints $V(U)$, they form a tree. We will construct $U$ in such a way it has no
finite end.

For $n \in \mathbb{N}$ set $E_n = [  \alpha n]$ be the integer part of
$\alpha N$, and set $e_n = E_n - E_{n - 1} \in \{ 0, 1 \}$ for $n \geqslant
1$. Given a vertex $x$ of $T$, remove one of the edges below it if $e_{| x| +
1} = 0$ and keep both of them if $e_{| x | + 1} = 1$. Let $U$ be the tree
obtained this way. 

  Let $n,p$ be in $\mathbb N$. Then, $E_{n + p} - E_p - n \alpha = O ( 1)$ as $n \rightarrow \infty$, uniformly
  in $p$:
\begin{equation}\label{lemma-model}
  E_{n+p} - E_p - n \alpha  =  [ ( n + p) \alpha] - ( n + p) \alpha - ( [
  p \alpha] - p \alpha)\in  [ - 1, 1] .
\end{equation}
Now, if $S_k(x)=\{y\succ x: |y-x|=k\}$,
$$
\log_2(\sharp S_k(x)/2^{\alpha k})= e_{|x|+1}+\ldots+e_{|x|+k}-\alpha k= E_{|x|+k}- E_{|x|}-\alpha k,
$$
hence, by \eqref{lemma-model},
$$
\frac14\leq \frac{\sharp S_k(x)}{2^{\alpha|x|}}\leq 1
$$
and the $\alpha$-regularity of $\partial U$ follows from Theorem \ref{ahl-regularity-trees}.

\vskip0.5cm

The proof of Theorem \ref{theomain} is similar, but it has to take the dishomogeneities of the tree into account.
We now prove Theorem \ref{theomain} for a general tree. We first need to homogenize the given tree $T$. Let $k>0$ be a fixed 
integer, then we construct a new tree $U$ which is easier to work with. This step could be avoided, but at the 
expenses of cleanness of exposition.

Consider the
portion $\Delta_{o, k}$ of the tree between $x_0 = o$ and $S_k ( x_0) = \{
x_{1, j} : | x_o - x_{1, j} | = k \}$ and replace it by geodesics of length $k$, each joining $x_o$ and some $x_{1, j}$, 
meeting at $x_o$ only. Repeat now
the same construction at each of the points $x_{1, j}$, considering the
portion of the tree between $x_{1, j}$ and $S_k ( x_{1, j})$ instead and iterate the process.  Let $U$ be the tree constructed 
this way. Let $S^T_{Nk}$ and $S^{U}_{Nk}$ be $S_{Nk}(x_0)$ for the tree $T$ and $U$ respectively. Then, we have
$$
S_{N k}^T = S_{N k}^U 
$$
for all $N \geqslant 1$. In fact, we have more. Let $\{x_{N, j_N}\}_N\subseteq S_{N k}$ be a sequence of points and let 
$T(x_{N,j_N})$ and $U(x_{N,j_N}) $ be the subtrees with root $x_{N,j_N}$ of $T$ and $U$ respectively. Then, 
$$
F : \bigcap_{N \geqslant 1} T ( x_{N, j_N}) \mapsto \bigcap_{N \geqslant 1}
   U ( x_{N, j_N}) , \ F:\partial T\to\partial U,
$$
is a bijection of $\partial T$ onto $\partial U$. Moreover, it is a bi-Lipschitz map.
\begin{lemma}\label{lip-tree}
  \begin{enumeratealpha}
    \item The map $F$ is a homeomorphism.
    
    \item Fix $0 < \delta \leqslant 1$. Then $F$ is bi-Lipschitz w.r.t. the
    metrics $\bar\rho^T$ on $\partial T$ and $\bar\rho^U$ on
    $\partial U$, with constants depending on $k\in\mathbb Z$.
  \end{enumeratealpha}
\end{lemma}
\begin{proof}
Assertion a) is implied by b). Let $\gamma$ be a geodesic joining $x, y
\in \partial T$. Then, $\gamma$ crosses a sequence of regions $\Delta_{x_j,k}^T$, entering and leaving in two points $y_{j - 1},
y_j \in \{ x_j \} \cup S_k^T ( x_j)$. The length of the geodesic
connecting the same vertices $y_{j - 1}, y_j$ in $\Delta_{x_j,k}^U $
is bounded above and below by positive multiplicative constants (depending on
$k$) times the corresponding quantity in $\Delta_{x_j,k}^T $. Hence,
$$
C_k^{- 1} \bar\rho^T (x, y) \leqslant
   \bar\rho^U ( F(x), F(y)) \leqslant C_k\bar\rho^T (
   x, y) 
$$
as we wished to show.
\end{proof}

\begin{figure}[H]
\begin{minipage}[b]{0.5\linewidth}
\centering

\begin{tikzpicture}[scale=0.9]
\draw[red] (0,0) --(-2,-1) ;
\draw (-2,-1) --(-3,-2);
\draw (-2,-1) --(-2,-2);
\draw[red] (-2,-1) --(-1,-2);

\draw[dotted] (-3,-2) --(-3,-5/2);
\draw[dotted] (-2,-2) --(-2,-5/2);

\draw[red] (-1,-2) --(-2,-7/2);
\draw (-1,-2) --(-0,-7/2);
\draw (-2,-7/2) --(-5/2,-9/2);
\draw[red] (-2,-7/2) --(-3/2,-9/2);

\draw (0,-7/2) --(-1/2,-9/2);
\draw (0,-7/2) --(1/2,-9/2);

\draw[dotted] (-1/2,-9/2) --(-1/2,-5);
\draw[dotted] (1/2,-9/2) --(1/2,-5);

\draw[dotted] (-5/2,-9/2) --(-5/2,-5);
\draw[dotted,red] (-3/2,-9/2) --(-3/2, -13/2);

\draw[fill] (0,0) circle [radius=0.05];
\draw[fill] (-2,-1) circle [radius=0.05];
\draw[fill] (-3,-2) circle [radius=0.05];
\draw[fill] (-2,-2) circle [radius=0.05];
\draw[fill] (-1,-2) circle [radius=0.05];
\draw[fill] (-2,-7/2) circle [radius=0.05];
\draw[fill] (0,-7/2) circle [radius=0.05];
\draw[fill] (-5/2,-9/2) circle[radius=0.05];
\draw[fill] (-3/2,-9/2) circle[radius=0.05];
\draw[fill] (-1/2,-9/2) circle[radius=0.05];
\draw[fill] (1/2,-9/2) circle[radius=0.05];
\draw[fill] (-3/2,-13/2) circle[radius=0.05];

\draw[dashed] (-7/2,-13/2) -- (-3/2,-13/2) --(7/2,-13/2);

\draw[dashed] (-7/2,-9/2) --(-5/2,-9/2) --(-3/2,-9/2) --(-1/2,-9/2) --(1/2,-9/2) --(7/2,-9/2);

\draw[dashed] (-7/2,-2) --(-2,-2) -- (-1,-2) --
 (1,-2) --(7/2,-2);
\node[above] at (0,0) {$o$};
\node[above] at (-5/2,-1/2) {\large $T$};
\node[above] at (-7/2,-13/2) {\footnotesize $\partial T$};
\node[above] at (-5/4,-13/2) {$x$};
\node[above] at (-7/2,-9/2) {\footnotesize $2k$};
\node[above] at (-7/2,-2) {\footnotesize $k$};
\draw (0,0) --(2,-1);
\draw (2,-1) --(1,-2);
\draw (2,-1) --(3,-2);

\draw[dotted] (1,-2) --(1,-5/2);
\draw[dotted] (3,-2) --(3,-5/2);

\draw[fill] (2,-1) circle [radius=0.05];
\draw[fill] (1,-2) circle [radius=0.05];
\draw[fill] (3,-2) circle [radius=0.05];
\end{tikzpicture}
\end{minipage}
\begin{minipage}[b]{0.5\linewidth}
\centering

\begin{tikzpicture}[scale=0.9]
\draw[red] (0,0) --(-1,-2) ;
\draw (0,0) --(-3,-2) ;
\draw (0,0) --(-2,-2) ;

\draw[dotted] (-3,-2) --(-3,-5/2);
\draw[dotted] (-2,-2) --(-2,-5/2);

\draw (-1,-2) --(-5/2,-9/2);
\draw[red] (-1,-2) --(-3/2,-9/2);
\draw (-1,-2) --(-1/2,-9/2);
\draw (-1,-2) --(1/2,-9/2);

\draw[dotted] (-1/2,-9/2) --(-1/2,-5);
\draw[dotted] (1/2,-9/2) --(1/2,-5);

\draw[dotted] (-5/2,-9/2) --(-5/2,-5);
\draw[dotted,red] (-3/2,-9/2) --(-3/2, -13/2);

\draw[fill] (0,0) circle [radius=0.05];
\draw[fill] (-1/2,-1) circle [radius=0.05];
\draw[fill] (-1,-1) circle [radius=0.05];
\draw[fill] (-3/2,-1) circle [radius=0.05];
\draw[fill] (-3,-2) circle [radius=0.05];
\draw[fill] (-2,-2) circle [radius=0.05];
\draw[fill] (-1,-2) circle [radius=0.05];
\draw[fill] (-5/4,-13/4) circle [radius=0.05];
\draw[fill] (-7/4,-13/4) circle [radius=0.05];
\draw[fill] (-3/4,-13/4) circle [radius=0.05];
\draw[fill] (-1/4,-13/4) circle [radius=0.05];
\draw[fill] (-5/2,-9/2) circle[radius=0.05];
\draw[fill] (-3/2,-9/2) circle[radius=0.05];
\draw[fill] (-1/2,-9/2) circle[radius=0.05];
\draw[fill] (1/2,-9/2) circle[radius=0.05];
\draw[fill] (-3/2,-13/2) circle[radius=0.05];

\draw[dashed] (-7/2,-13/2) -- (-3/2,-13/2) --(7/2,-13/2);

\draw[dashed] (-7/2,-9/2) --(-5/2,-9/2) --(-3/2,-9/2) --(-1/2,-9/2) --(1/2,-9/2) --(7/2,-9/2);

\draw[dashed] (-7/2,-2) --(-2,-2) -- (-1,-2) --
 (1,-2) --(7/2,-2);
\node[above] at (0,0) {$o$};
\node[above] at (-5/2,-1/2) {\large $U$};
\node[above] at (7/2,-13/2) {\footnotesize $\partial U$};
\node[above] at (-1,-13/2) {$F(x)$};
\node[above] at (7/2,-9/2) {\footnotesize $2k$};
\node[above] at (7/2,-2) {\footnotesize $k$};

\draw (0,0) --(1,-2);
\draw (0,0) --(3,-2);

\draw[dotted] (1,-2) --(1,-5/2);
\draw[dotted] (3,-2) --(3,-5/2);

\draw[fill] (1,-2) circle [radius=0.05];
\draw[fill] (3,-2) circle [radius=0.05];
\draw[fill] (1/2,-1) circle [radius=0.05];
\draw[fill] (3/2,-1) circle [radius=0.05];
\end{tikzpicture}
\end{minipage}
\caption{From the tree $T$ to the tree $U$.}
\end{figure}
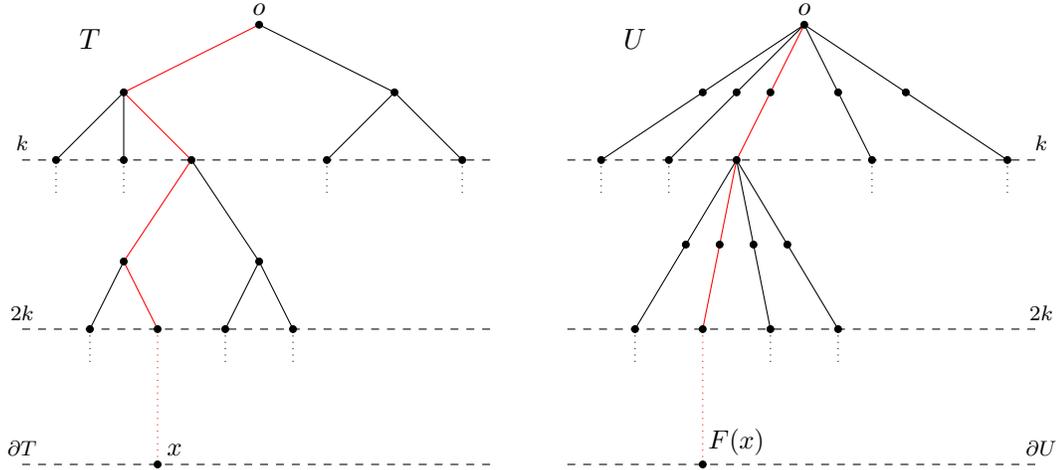
Since Ahlfors regularity is a metric property, preserved by bi-Lipschitz maps, we might look for our Ahlfors regular 
subspace on $\partial U$ rather than on $\partial T$.

We now prove the dyadic version of Theorem \ref{theomain}.
\begin{thm}\label{main-tree}
Let $(\partial T,\bar\rho)$ be a Ahlfors $Q$-regular metric space and let $0<\alpha<Q$. Then, there exists 
$Y\subseteq \partial T$ such that $(Y,\bar\rho)$ is $\alpha$-regular. 
\end{thm}

\begin{proof} Given the data $0<\alpha<Q$, let $\epsilon>0$ such that $\alpha< Q-\epsilon$ is rational and let $k>0$ integer such
that (i) $(Q-\epsilon)k$ is integer, and (ii) $ \frac{\sharp S^T_k(x)}{2^{kQ}}\geq 2^{-\varepsilon k}$ (this requirement can be 
fulfilled if $k$ is large enough, by Theorem \ref{ahl-regularity-trees}). 

 For such a $k$, let $U$ be the tree in Lemma \ref{lip-tree}, 
 whose boundary is biLipschitz equivalent to the boundary of $T$. 
 Clearly, the metric space $(\partial U,\bar\rho)$ is Ahlfors $Q$-regular as well. We now construct a 
 subtree of $U$ as follows. Given the root $o\in U$, we select exactly $2^{(Q-\varepsilon)k}$ vertices in the set 
 $S^U_{k}(o)=\{y\succ x: |y-o|=k\}$. Let $\{x_{1,j}\}_{j=1}^{2^{(q-\varepsilon)k}}$ be these points. Then, we select 
 $2^{(Q-\varepsilon)k}$ vertices in the set $S^{U}_{k}(x_{1,j})$ for any $j=1,\ldots 2^{(Q-\varepsilon)k}$ and we iterate the 
 process. We denote by $V$ the subtree of $U$ constructed in such a way. The tree $V$ is very close to be a homogeneous tree, 
 it is a \emph{periodic} tree. The \emph{level} $n$ of the tree $V$ is the set of points $\{x\in V: |x|= nk\}$
    
We now proceed to the construction of an $\alpha$-regular subset $\partial W$ of
$\partial V$. The desired $\alpha$-regular subspace of $\partial T$ will be $F^{- 1} ( \partial W)$ where $F$ is the map of 
Lemma \ref{lip-tree}.

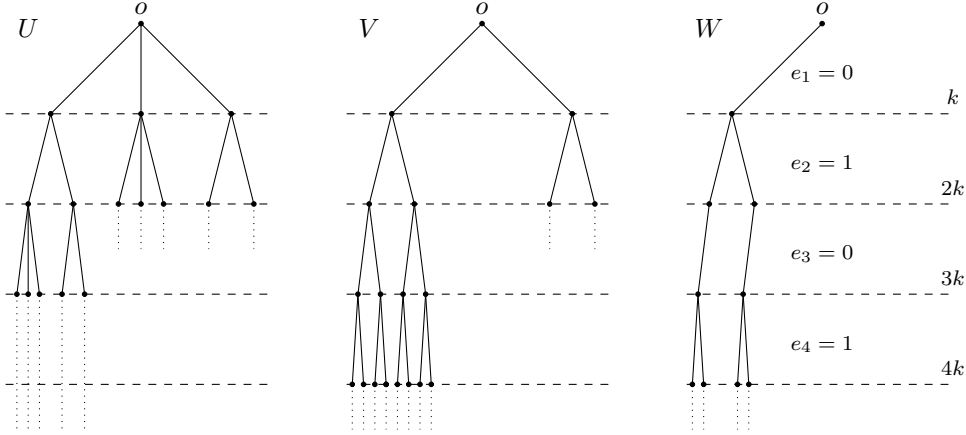
\begin{figure}
\begin{minipage}[t]{0.3\linewidth}
\centering
\begin{tikzpicture}[scale=0.6]
\draw (0,0) --(-2,-2);
\draw (0,0) --(0,-2);
\draw (0,0) --(2,-2);

\draw (-2,-2) --(-5/2,-4);
\draw (-2,-2) --(-3/2,-4);

\draw (0,-2) --(-1/2,-4);
\draw (0,-2) --(0,-4);
\draw (0,-2) --(1/2,-4);

\draw (2,-2) --(3/2,-4);
\draw (2,-2) --(5/2,-4);

\draw[dotted] (-1/2,-4) --(-1/2,-5);
\draw[dotted] (0,-4) --(0,-5);
\draw[dotted] (1/2,-4) --(1/2,-5);

\draw[dotted] (3/2,-4) --(3/2,-5);
\draw[dotted] (5/2,-4) --(5/2,-5);
 
\draw (-5/2,-4) --(-11/4,-6);
\draw (-5/2,-4) --(-5/2,-6);
\draw (-5/2,-4) --(-9/4,-6);

\draw (-3/2,-4) --(-7/4,-6);
\draw (-3/2,-4) --(-5/4,-6);

\draw[dotted] (-11/4,-6) --(-11/4,-9);
\draw[dotted] (-5/2,-6) --(-5/2,-9);
\draw[dotted] (-9/4,-6) --(-9/4,-9);
\draw[dotted] (-7/4,-6) --(-7/4,-9);
\draw[dotted] (-5/4,-6) --(-5/4,-9);

\draw[fill] (0,0) circle [radius=0.05];
\draw[fill] (-2,-2) circle [radius=0.05];
\draw[fill] (0,-2) circle [radius=0.05];
\draw[fill] (2,-2) circle [radius=0.05];
\draw[fill] (-1/2,-4) circle [radius=0.05];
\draw[fill] (0,-4) circle [radius=0.05];
\draw[fill] (1/2,-4) circle [radius=0.05];
\draw[fill] (-5/2,-4) circle [radius=0.05];
\draw[fill] (-3/2,-4) circle [radius=0.05];
\draw[fill] (3/2,-4) circle [radius=0.05];
\draw[fill] (5/2,-4) circle [radius=0.05];
\draw[fill] (-11/4,-6) circle [radius=0.05];
\draw[fill] (-9/4,-6) circle [radius=0.05];
\draw[fill] (-5/2,-6) circle [radius=0.05];
\draw[fill] (-7/4,-6) circle [radius=0.05];
\draw[fill] (-5/4,-6) circle [radius=0.05];

\draw[dashed] (-3,-8) --(2.9,-8);
\draw[dashed] (-3,-6) --(2.9,-6);
\draw[dashed] (-3,-4) --(2.9,-4);
\draw[dashed] (-3,-2) --(2.9,-2);
\node[above] at (0,0) {$o$};
\node[above] at (-5/2,-1/2) { $U$};

\end{tikzpicture}
\end{minipage}
\begin{minipage}[t]{0.3\linewidth}
\centering
\begin{tikzpicture}[scale=0.6]
\draw (0,0) --(-2,-2);
\draw (0,0) --(2,-2);

\draw (-2,-2) --(-5/2,-4);
\draw (-2,-2) --(-3/2,-4);


\draw (2,-2) --(3/2,-4);
\draw (2,-2) --(5/2,-4);


\draw[dotted] (3/2,-4) --(3/2,-5);
\draw[dotted] (5/2,-4) --(5/2,-5);
 
\draw (-5/2,-4) --(-11/4,-6);
\draw (-5/2,-4) --(-9/4,-6);

\draw (-3/2,-4) --(-7/4,-6);
\draw (-3/2,-4) --(-5/4,-6);

\draw (-11/4,-6) --(-23/8,-8);
\draw (-11/4,-6) --(-21/8,-8);
\draw (-9/4,-6) --(-17/8,-8);
\draw (-9/4,-6) --(-19/8,-8);

\draw (-7/4,-6) --(-15/8,-8);
\draw (-7/4,-6) --(-13/8,-8);
\draw (-5/4,-6) --(-11/8,-8);
\draw (-5/4,-6) --(-9/8,-8);

\draw[dotted] (-9/8,-8) --(-9/8,-9);
\draw[dotted] (-11/8,-8) --(-11/8,-9);
\draw[dotted] (-13/8,-8) --(-13/8,-9);
\draw[dotted] (-15/8,-8) --(-15/8,-9);
\draw[dotted] (-19/8,-8) --(-19/8,-9);
\draw[dotted] (-21/8,-8) --(-21/8,-9);
\draw[dotted] (-17/8,-8) --(-17/8,-9);
\draw[dotted] (-23/8,-8) --(-23/8,-9);

\draw[fill] (0,0) circle [radius=0.05];
\draw[fill] (-2,-2) circle [radius=0.05];
\draw[fill] (2,-2) circle [radius=0.05];
\draw[fill] (-5/2,-4) circle [radius=0.05];
\draw[fill] (-3/2,-4) circle [radius=0.05];
\draw[fill] (3/2,-4) circle [radius=0.05];
\draw[fill] (5/2,-4) circle [radius=0.05];
\draw[fill] (-11/4,-6) circle [radius=0.05];
\draw[fill] (-9/4,-6) circle [radius=0.05];
\draw[fill] (-7/4,-6) circle [radius=0.05];
\draw[fill] (-5/4,-6) circle [radius=0.05];
\draw[fill] (-23/8,-8) circle [radius=0.05];
\draw[fill] (-21/8,-8) circle [radius=0.05];
\draw[fill] (-19/8,-8) circle [radius=0.05];
\draw[fill] (-17/8,-8) circle [radius=0.05];
\draw[fill] (-17/8,-8) circle [radius=0.05];
\draw[fill] (-15/8,-8) circle [radius=0.05];
\draw[fill] (-13/8,-8) circle [radius=0.05];
\draw[fill] (-11/8,-8) circle [radius=0.05];
\draw[fill] (-9/8,-8) circle [radius=0.05];

\draw[dashed] (-3,-8) --(2.9,-8);
\draw[dashed] (-3,-6) --(2.9,-6);
\draw[dashed] (-3,-4) --(2.9,-4);
\draw[dashed] (-3,-2) --(2.9,-2);
\node[above] at (0,0) {$o$};
\node[above] at (-5/2,-1/2) { $V$};

\end{tikzpicture}
\end{minipage}
\begin{minipage}[t]{0.3\linewidth}
\centering
\begin{tikzpicture}[scale=0.6]
\draw (0,0) --(-2,-2);

\draw (-2,-2) --(-5/2,-4);
\draw (-2,-2) --(-3/2,-4);




 
\draw (-5/2,-4) --(-11/4,-6);

 \draw (-3/2,-4) --(-7/4,-6);

\draw (-11/4,-6) --(-23/8,-8);
\draw (-11/4,-6) --(-21/8,-8);

\draw (-7/4,-6) --(-15/8,-8);
\draw (-7/4,-6) --(-13/8,-8);

\draw[dotted] (-13/8,-8) --(-13/8,-9);
\draw[dotted] (-15/8,-8) --(-15/8,-9);
\draw[dotted] (-21/8,-8) --(-21/8,-9);
\draw[dotted] (-23/8,-8) --(-23/8,-9);

\draw[fill] (0,0) circle [radius=0.05];
\draw[fill] (-2,-2) circle [radius=0.05];
\draw[fill] (-5/2,-4) circle [radius=0.05];
\draw[fill] (-3/2,-4) circle [radius=0.05];
\draw[fill] (-11/4,-6) circle [radius=0.05];
\draw[fill] (-7/4,-6) circle [radius=0.05];
\draw[fill] (-23/8,-8) circle [radius=0.05];
\draw[fill] (-21/8,-8) circle [radius=0.05];
\draw[fill] (-15/8,-8) circle [radius=0.05];
\draw[fill] (-13/8,-8) circle [radius=0.05];

\draw[dashed] (-3,-8) --(2.9,-8);
\draw[dashed] (-3,-6) --(2.9,-6);
\draw[dashed] (-3,-4) --(2.9,-4);
\draw[dashed] (-3,-2) --(2.9,-2);
\node[above] at (0,0) {$o$};
\node[above] at (-5/2,-1/2) { $W$};
\node[above] at (2.9,-8) {\footnotesize $4k$};
\node[above] at (2.9,-6) {\footnotesize $3k$};
\node[above] at (2.9,-4) {\footnotesize $2k$};
\node[above] at (2.9,-2) {\footnotesize $k$};
\node[above] at (0,-3/2) {\footnotesize $e_1=0$};
\node[above] at (0,-7/2) {\footnotesize $e_2=1$};
\node[above] at (0,-11/2) {\footnotesize $e_3=0$};
\node[above] at (0,-15/2) {\footnotesize $e_4=1$};


\end{tikzpicture}
\end{minipage}
\caption{The trees $U$, $V$ and $W$.}
\end{figure}

Given $0<\alpha<Q$, 
set $\lambda=\frac{\alpha}{Q-\epsilon}\in(0,1)$.
Let $E_N=[\lambda N]$, and let
$\{ e_n=E_n-E_{n-1} : n \geqslant 1 \}$, so that
 $e_n \in \{ 0, 1 \}$. Construct a subtree $W$ of $V$ accordingly to the following procedure. 
For each vertex $x$ at the $n^{th}$ level of $V$: 
if $e_n=1$, we keep all the edges leaving $x$ in direction of the $(n+1)^{th}$ level; 
if $e_n=0$ we dismiss all the edges leaving $x$ in direction of the $(n+1)^{th}$ level  but one. 
We now prove that the tree $W$ we constructed is $\alpha$-regular exploiting Theorem \ref{ahl-regularity-trees}.

Recall that $E_n = e_1 + \ldots + e_n$, $n\in\mathbb N$. 
Fix now $x$ in $W$, $(p-1)k<|x|\le pk$, and let $l\ge pk-|x|$,
$$
l=[pk-|x|]+nk+m, \text{ with }0\le m<k.
$$
Then,
\begin{eqnarray*}
 \sharp(S^W_l(x))&=&\sharp(S^W_{nk}(x_1)) \text{ where $x_1$ is the point in $W$ below $x$ such that $|x_1|=pk$}\crcr
 &=&2^{(Q-\epsilon)k(E_{n+p}-E_p)} \text{ by construction of $W$}\crcr
 &\approx&2^{(Q-\epsilon)kn\lambda} \text{ by \eqref{lemma-model}},\crcr
 &=&2^{\alpha kn}\crcr
 &\approx&2^{\alpha l},
\end{eqnarray*}
with constants which depend on $k$. For $0<l<pk-|x|$, $\sharp(S^W_l(x))=1\approx 2^{\alpha l}$, 
hence the condition in Theorem \ref{ahl-regularity-trees} is satisfied, showing that $\partial W$ is $\alpha$-regular,
and $F^{-1}(\partial W)$ is the $\alpha$-regular subspace of $\partial T$ we were looking for.
\end{proof}

The constant $C>0$ which satisfies 
$$
 2^{\alpha\ell}/C\leq \sharp S^W_\ell(x)\leq C 2^{\alpha\ell}
$$
for all $x\in W$ and $\ell\in\mathbb N$ depends on the integer $k$ used to construct the tree 
$U$ in the proof of Theorem \ref{main-tree}. It would be interesting to have better quantitative information
on the dependence of $C$ on the data.

  \section{From trees to metric spaces.}\label{dyadization}
  
  We now proceed to deduce Theorem \ref{theomain} in the general context of Ahlfors regular spaces, from the special tree case. 
  A key role is played by M. Christ's dyadic decomposition of metric spaces, which we now recall.

Suppose the environment space $(X,\rho)$ is an Ahlfors regular metric measure space having dimension $Q$. 
It easily follows from $Q$-regularity that
\begin{equation}\label{diametro}
 \text{diam}(B(x,r))\approx r.
\end{equation}
Michael Christ gave, in the more general context of metric spaces of homogeneous type, 
a dyadic decomposition of $X$. We summarize here Christ's result \cite{Christ}
in the context of Ahlfors regular spaces.
\begin{lemma}\label{lemmachrist}
Suppose $(X,\rho)$ is $Q$-regular and that $\diam(X)=1$. 
There exists a collection $\{\icirc^k_a, a\in A_k, k\in\mathbb Z\}$ ($A_k$ being a finite set of indices) 
of open subsets of $X$ and $0<\delta<1$ and $c_1,c_2>0$ 
such that: 
\begin{enumerate}[(i)]
 \item $\HHH^Q\left(X\setminus\bigcup_{a \in A_k} \icirc_a^k\right)=0$ holds for all $k\geq 0$;
 \item if $l\geq k$ and $a\in A_k$, then for all $b\in A_l$ either $\icirc_b^l\subseteq\icirc_a^k$ or $\icirc_a^k\cap\icirc_b^l=\emptyset$;
 \item if $l> k$ and $b\in A_l$, then there exists a unique $a\in A_k$ such that $\icirc_b^l\subseteq\icirc_a^k$;
 \item  $\diam(\icirc_a^k)\le c_1\delta^k$ for some constant $c_1$ depending on the metric space;
 \item for all $(k,a)$, there are $z_a^k\in \icirc^k_a$ such that $B(z_a^k,c_2\delta^k)\subseteq\icirc_a^k$.
\end{enumerate}
\end{lemma}
From inspection of the proof, it is clear that any choice of $\delta, c_2$ such that $\delta+c_2\le\frac{1}{4}$ will work in 
general.

We can assume that $\sharp A_0=1$. 
From now on we set $I_a^k=\overline{\icirc_a^k}$ and let $T=\cup_{k\geq 0}\left\{I_a^k:\ a\in A_k\right\}$. 
Sometimes we simply write $I_a^k=I\in A_k$.
The set $T$ has a tree structure, where there is a undirected edge $((k,a),(k+1,b))=(\alpha,\beta)$ if 
$I_a^k\subseteq I_b^{k+1}$.  

Given $\delta\in(0,1)$, we can introduce a distance $\rho_\delta$ on $T$ by assigning the weight $\delta^k$ to the edge 
$((k,a),(k+1,b))$, measuring lengths of paths according to the weight,
and defining the distance between two vertices to be the shortest length of a path joining them, 
as we did in Section \ref{alberi},
 where we had set $\delta=\frac{1}{2}$. The change of parameter has little consequences. 
 The Cauchy sequences are the same, and so are the completion $\overline{X}$ and the boundary $\partial T$. 
 The only difference is 
 that $(\partial T,\rho_\delta)$ and $(\partial T,\rho_{1/2})$ are ``snowflake versions'' of each other: there is $C>0$
 such that, with $\rho_{1/2}=\bar\rho$,
 $$
 \rho_\delta(x,y)=C\rho_{1/2}(x,y)^{\log_2(1/\delta)},
 $$
 hence,
\begin{eqnarray*}
 B_{\rho_\delta}(x,r)&=&B_{\rho_{1/2}}(x,C r^{\log_2(1/\delta)}),\crcr
 \text{diam}_{\rho_\delta}(U)&=&\text{diam}_{\rho_{1/2}}(U)^{\log_2(1/\delta)}, \text{ and }\crcr
  \text{dim}_{\rho_{\delta}}(\partial T)&=&=\log_2(1/\delta)\text{dim}_{\rho_{1/2}}(\partial T),
\end{eqnarray*}
 statements which are easily verified using the ultrametric property of the distances $\rho_\delta$.

We consider the map $\Lambda:\partial T\to X$, associating to 
 $x=\{(k,a)\}_{k\geq0}\in\partial T$ the element 
$\Lambda(x)=\cap_{k\geq0}I^k_a\in X$. As proved in \cite{ARSW}, the map $\Lambda$ is a surjective map with some 
nice properties. In particular,
\begin{lemma}\label{whatWeNeed}
\begin{enumerate}
 \item[(a)] $\Lambda:(\partial T,\rho_\delta)\to (X,\rho)$ is Lipschitz.
 \item[(b)] $(\partial T,\rho_\delta)$ is $Q$-regular if $(X,\rho)$ is $Q$-regular..
 \item[(c)] For each $k\ge0$ there is a bijection $G_k:S_k(o)\to A_k$ so that, if $a=G_k(x)$, then $I_a^k=\Lambda(\partial T(x))$. 
In this case, we write $\partial T(x)=\widetilde{I}_a^k$.
\end{enumerate}
\end{lemma}
A key fact for us is that the map $\Lambda$ is \emph{regular} in the sense of G. David and S. Semmes \cite{DS-mon, DS}.
\begin{defn}
 Let $(X, \rho )$ and $(Y, \nu)$ be two metric spaces.
A mapping $f : X\to Y$ is said to be \emph{regular} if it is Lipschitz, and if there is a constant $C > 0$ such that for every ball $B$ in $Y$
it is possible to cover $f^{-1}(B)$ by at most $C$ balls in $X$ with radius equal to $C \cdot radius(B)$.
\end{defn}
Once we prove that $\Lambda$ is a regular map, the main theorem will follow from the following lemma from 
\cite[Lemma 3.10]{DS}.
\begin{lemma}\label{DS}
Let $(X, \rho)$ and $(Y, \nu)$ be metric spaces, and let
$f : X \to Y$ be a regular map. If $X$ is Ahlfors-regular of  dimension $Q$, then $f(X)$ is Ahlfors-
regular of dimension $Q$ as well.
\end{lemma}
It suffices to prove the following.
\begin{prop}\label{lambda-regular1}
\begin{enumerate}
 \item[(a)] The map $\Lambda:\partial T\to X$ is regular.
 \item[(b)]  Let $Y\subset\partial T$ a $\alpha$-regular subspace of $\partial T$. 
 Then, the map $\Lambda:Y\to \Lambda(Y)$ is regular. 
\end{enumerate}
\end{prop}
\begin{proof} (a)
We first prove that there exists a constant $c_3>0$ such that, given any ball $B:=B(x,r)\subseteq X$ with 
$\delta^{k+1}<r\le\delta^k$, 
$B\subseteq\cup_{j=1}^{c_3}I_{j}$, where the $I_j$'s are Christ's cubes with labels in $A_k$. 
Consider  the set 
$\III=\{I_j: I_j\cap B\neq\emptyset, I_j\in A_k\}$ and set $N:=\sharp \III$. 
Then, there exist a constant $c>0$ and a ball $B'$ of radius $c\delta^k$ such that 
$\cup_{I_j\in \III}I_j\subseteq B'$. Now, from $Q$-regularity and $(v)$ 
in Lemma \ref{lemmachrist}, we get that 
for any $I_j\in \III$, $\mathcal H^{Q}(I_j)\gtrsim \delta^{Qk}$. 
Therefore, again by $Q$-regularity,
\begin{align*}
 N \delta^{Qk}\lesssim\sum_{I_j\in \III}\mathcal H^{Q}(I_j)\lesssim \mathcal H^{Q}(B')\lesssim  \delta^{Qk}.
\end{align*}
Hence, there exists a constant $c_3>0$ such that $N< c_3<+\infty$.

For each dyadic cube $I$ in $X$, denote by $\widetilde I\subseteq \partial T$ the corresponding dyadic cube (and metric 
cl-open ball) in $\partial T$.
Now we prove that given any dyadic cube  $I\subseteq X$ with $I\in A_k$, the set 
$\widetilde \III:=\{\widetilde I_j\in A_k:  \widetilde {I_j}\subseteq\partial T \text{ such that } 
\Lambda(\widetilde I_j)\cap I\neq\emptyset\}$ has bounded cardinality. 
 Set $I_j:=\Lambda(\widetilde I_j)$ and consider $\cup_{j}I_j$. All the cubes $I_j$ intersect 
the given cube $I\subset X$ and, reasoning as in the first part of the prove, we obtain that the number of cubes $I_j$ 
has to be bounded by an absolute constant independent of the given cube $I$. Therefore, the number of cubes $\widetilde I_j$ 
has to be bounded as well by an absolute constant, say $c_4>0$. To recap, we proved that any ball $B\subseteq X$ can be 
covered by at most $c_3$ dyadic cubes whose diameter is comparable to the radius of the ball $B$ and to any dyadic cube in $X$ 
we associated, by means of $\Lambda$, at most $c_4$ dyadic cubes in $\partial T$ with the same diameter. Recalling that dyadic 
cubes and balls coincides in $\partial T$, we conclude that $\Lambda^{-1}(B)$ can be covered by at most $c_3c_4$ balls in 
$\partial T$ with radius comparable to the one of $B$. This concludes the proof.

(b)
Let $x\in \Lambda(Y)$ and consider the ball $B\subseteq\Lambda(Y)$ given by $B:=B(x,r)\cap \Lambda(Y)$, where $B(x,r)$ is a
 ball in $X$. Then, we want to cover $\Lambda^{-1}(B(x,r)\cap \Lambda(Y))$ with a bounded number of balls in $Y$ whose radius is
 comparable to the radius of $B$. By (a), $\Lambda^{-1}(B(x,r))$ can be covered by a finite number of balls $\widetilde{I}_j$ in $\partial T$, having radius comparable to $r$. Discard those which do not intersect $Y$ and, for each $\emptyset\ne\widetilde{I}_j\cap Y\ni y_j$,
we can take $y_j$ to be its center, by the ultrametric property. Such sets
$\widetilde{I}_j\cap Y$ form a covering of $\Lambda^{-1}(B(x,r)\cap \Lambda(Y))$ by
boundedly many balls in $Y$, having radius comparable to $r$.  and we conclude that $\Lambda: Y\to \Lambda({Y})$ is a regular map as wished.
 \end{proof}

To finish the proof, consider 
$$
0<\alpha<Q=\text{dim}_\rho(X)=\text{dim}_{\rho_\delta}(\partial T)=\log_2(1/\delta)\text{dim}_{\rho_{1/2}}(\partial T).
$$
By Theorem \ref{ahl-regularity-trees}, there is an Ahlfors regular subspace $Y$ of $\partial T$ with
$$
\frac{\alpha}{\log_2(1/\delta)}=\text{dim}_{\rho_{1/2}}(Y)=\frac{\text{dim}_{\rho_\delta(Y)}}{\log_2(1/\delta)}.
$$
By Proposition \ref{lambda-regular1} (b), $\Lambda(Y)$ is an $\alpha$-regular subspace of $X$, as stated in 
Theorem \ref{theomain}.

\vskip0.5cm

It is a pleasure to thank Yuval Peres for directing us to the reference \cite{How95}, and Pekka Koskela who made us aware of \cite{JJKRRS} which was not cited in an earlier version of our paper.

\bibliography{ahlfors-bib}

\newcommand{\etalchar}[1]{$^{#1}$}
\providecommand{\bysame}{\leavevmode\hbox to3em{\hrulefill}\thinspace}
\providecommand{\MR}{\relax\ifhmode\unskip\space\fi MR }
\providecommand{\MRhref}[2]{%
  \href{http://www.ams.org/mathscinet-getitem?mr=#1}{#2}
}
\providecommand{\href}[2]{#2}
\begin{thebibliography}{ARSW14}

\bibitem[ARSW14]{ARSW}
N.~Arcozzi, R.~Rochberg, E.~T. Sawyer, and B.~D. Wick, \emph{Potential theory
  on trees, graphs and {A}hlfors-regular metric spaces}, Potential Anal.
  \textbf{41} (2014), no.~2, 317--366.

\bibitem[Chr90]{Christ}
M.~Christ, \emph{A {$T(b)$} theorem with remarks on analytic capacity and the
  {C}auchy integral}, Colloq. Math. \textbf{60/61} (1990), no.~2, 601--628.

\bibitem[DS97]{DS-mon}
G.~David and S.~Semmes, \emph{Fractured fractals and broken dreams:
  self-similar geometry through metric and measure}, Oxford Lecture Series in
  Mathematics and its Applications, vol.~7, The Clarendon Press, Oxford
  University Press, New York, 1997.

\bibitem[DS00]{DS}
\bysame, \emph{Regular mappings between dimensions}, Publ. Mat. \textbf{44}
  (2000), no.~2, 369--417.

\bibitem[Fal86]{Fal}
K.J. Falconer, \emph{The geometry of fractal sets}, Cambridge Tracts in
  Mathematics, vol.~85, Cambridge University Press, Cambridge, 1986.

\bibitem[HK98]{HK}
J.~Heinonen and P.~Koskela, \emph{Quasiconformal maps in metric spaces with
  controlled geometry}, Acta Math. \textbf{181} (1998), no.~1, 1--61.

\bibitem[How95]{How95}
J.~D. Howroyd, \emph{On dimension and on the existence of sets of finite
  positive {H}ausdorff measure}, Proc. London Math. Soc. (1995).

\bibitem[JJK{\etalchar{+}}10]{JJKRRS}
E.~J{\"a}rvenp{\"a}{\"a}, M.~J{\"a}rvenp{\"a}{\"a}, A.~K{\"a}enm{\"a}ki,
  T.~Rajala, S.~Rogovin, and V.~Suomala, \emph{Packing dimension and {A}hlfors
  regularity of porous sets in metric spaces}, Mathematische Zeitschrift
  \textbf{266} (2010), no.~1, 83--105.

\bibitem[Rog70]{Rogers}
C.~A. Rogers, \emph{Hausdorff measures}, Cambridge University Press, London-New
  York, 1970.

\end{thebibliography}
\bibliographystyle{amsalpha}  

\end{document}